\newif\ifArxivVersion
\newif\ifShowProofs
\newif\ifUnshortenedProofs

\ArxivVersiontrue
\ShowProofstrue
\UnshortenedProofstrue

\documentclass[letterpaper, 10 pt, conference]{ieeeconf}  

\IEEEoverridecommandlockouts                              
\usepackage{cite}
\usepackage{textcomp}

\newif\ifshowproofs
\showproofstrue

\usepackage{cite}
\usepackage{hyperref}

\usepackage{amsmath}
\usepackage{amssymb}
\usepackage{mathtools}
\usepackage{bbm}
\usepackage{graphicx}
\usepackage{empheq}
\usepackage{tikz}
\usepackage{algpseudocode}
\usepackage{algorithm2e}
\usepackage{subcaption}

\SetKwComment{Comment}{/* }{ */}
\RestyleAlgo{ruled}
\usepackage{pgfplots}
\usepackage[inkscapeformat=png]{svg}

\newtheorem{theorem}{Theorem}[section]
\newtheorem{lemma}[theorem]{Lemma}
\newtheorem{proposition}[theorem]{Proposition}

\newtheorem{definition}[theorem]{Definition}

\newtheorem{assumption}{Assumption}

\usepackage{mathtools}
\usepackage{fancyhdr}
\usepackage{amsfonts}
\usepackage{amssymb}
\usepackage[T1]{fontenc}
\usepackage{siunitx}
\usepackage{multirow}
\usepackage{placeins}
\usepackage{booktabs}
\usepackage{rotating}
\usepackage{array}
\usepackage{dsfont}
\usepackage{bbold}
\usepackage{adjustbox}
\usetikzlibrary{calc}
\usepackage{xcolor}
\definecolor{niceGreen}{rgb}{0.1, 0.625, 0.1}

\DeclareMathOperator*{\argmin}{arg\,min}

\usepackage[capitalise]{cleveref}

\newcounter{Cequ}

  \usepackage{tikz}
\usetikzlibrary{automata, positioning, arrows}
\tikzset{
->, 
node distance=2.2cm, 
every state/.style={thick, fill=gray!10}, 
initial text=$ $, 
}

\tikzset{every loop/.style={min distance=1mm,in=75,out=105,looseness=5}}

\usetikzlibrary{arrows.meta}
\tikzset{>={Latex[width=2mm,length=2mm]}}

\tikzset{every path/.style={line width=0.3 mm}}
\usepackage{multirow}

\title{\LARGE \bf
Joint Chance Constrained Optimal Control via Linear Programming
}

\author{Niklas Schmid, Marta Fochesato, Tobias Sutter, John Lygeros
\thanks{N. Schmid, M. Fochesato, and J. Lygeros are with the Automatic Control Laboratory, ETH Zürich, 8092 Zürich, Switzerland \{\tt\small nschmid,mfochesato,jlygeros\}@ethz.ch}%
\thanks{T.~Sutter is with the Department of Computer Science, University of Konstanz, 78457 Konstanz, Germany
        {\tt\small \ tobias.sutter@uni-konstanz.de}}%
        \thanks{Work supported by the European Research Council under the Horizon 2020 Advanced under Grant 787845 (OCAL).}%
}

\usepackage[printonlyused]{acronym}
\usepackage{xcolor}
\usepackage{amsmath}
\usepackage{amssymb}
\usepackage{mathtools}
\usepackage{subcaption}
\usepackage{xcolor}
\usepackage{hyperref}
\hypersetup{colorlinks=true, linkcolor=black}
\usepackage{mathrsfs}
\usepackage{cite}
\usepackage{adjustbox}
\usepackage{float}

\DeclarePairedDelimiterX{\inp}[2]{\langle}{\rangle}{#1, #2}

\newcommand\mydots{\hbox to 1em{.\hss.\hss.}}

\begin{document}

\maketitle
\thispagestyle{empty}
\pagestyle{empty}

\begin{abstract}
We establish a linear programming formulation for the solution of joint chance constrained optimal control problems over finite time horizons. The joint chance constraint may represent an invariance, reachability or reach-avoid specification that the trajectory must satisfy with a predefined probability. For finite state and action spaces, the solution is exact and our method computationally superior to approaches in the literature. For continuous state or action spaces, our linear programming formulation enables basis function approximations. 
\end{abstract}


\section{Introduction}
\label{sec:introduction}
Most system dynamics are non-deterministic, either naturally or due to unobservable modes, e.g., variability of large scale processes \cite{pitchford2007uncertainty, fochesato2022data} or robotic path planning in uncertain environments \cite{schmid2022real, blackmore2011chance}. While control problems often feature a ``real" control cost objective, e.g., energy consumption or negative financial profit, additional objectives often need to be included, such as reaching or remaining in a set of states, e.g., for safety. Such behaviours can be imposed in optimal control problems by manually shaping the cost associated to individual states until the policy behaves as desired. However, without any instructive procedure it is generally hard to achieve clear trade-offs between safety guarantees and the control cost \cite{Ono_2}. 

The task of remaining in a safe set is called \textit{invariance}, that of reaching a target set \textit{reachability}, while \textit{reach-avoidance} combines both objectives. The goal of finding a policy that yields an optimal trade-off between the control cost incurred by the system and the probability of achieving one of the above specifications is described by the so called joint-chance constrained optimal control problem (cf. \cite{schmid2023Jccoc}). While this class of problems is very popular, existing solution approaches are computationally challenging due to the non-Markovian structure imposed by these constraints (cf. \cite{Wang}). To overcome these computational issues, the state-of-the-art is to utilize approximations based on conservative bounds \cite{Ono_2, paulson, patil_1, wang2020non, nemirovski2007convex, hokayem2013chance, blackmore2011chance} or approximate the stochastic elements via sampling \cite{thorpe_1
, huang2021risk}.

A recently evolving idea is to break the non-Markovian structure using a state-augmentation, resulting in a conventional Markov Decision Process (MDP) \cite{schmid2023Jccoc, haesaert2021formal, hahn2019interval}. In \cite{hahn2019interval}, a value iteration algorithm is formulated to solve the full Pareto front for interval MPDs with multiple probabilistic temporal logic objectives. In \cite{haesaert2021formal}, a linear program (LP) on occupation measures for interval MPDs is formulated to find trade-offs between multiple probabilistic temporal logic objectives. However, the proposed solutions feature only interval MDPs with finite states and actions and the potential suboptimality of the computed policies remains unaddressed. In \cite{schmid2023Jccoc}, the joint chance constrained problem is approached via a bilevel optimization on its Lagrange dual, requiring to iteratively solve Dynamic Programming (DP) recursions. A suboptimality bound to the global optimum is presented, which decreases exponentially with the number of DP iterations. 
Unfortunately, the suboptimality bound in \cite{schmid2023Jccoc} convergences to an optimal policy only in the limit, while it is desirable to obtain a provably optimal policy in finitely many computations. This paper closes this gap. 

\textbf{Contribution.} As our main contribution we show that the bilevel problem in \cite{schmid2023Jccoc} can be solved via a sequence of two LPs. We prove that the obtained policy is an optimal solution to the joint chance constrained optimal control problem. To extend the applicability of our results we provide instructions on how to extend the safety specification in \cite{schmid2023Jccoc} from invariance to reach-avoid specifications. Faced with large or continuous state spaces, the LP framework allows for basis function approximations with generally superior scalability compared to gridding~\cite{de2003linear}.

\textbf{Notation.} We denote by $[N]$ the set $\{0,1,\dots,N\},N\in\mathbb{N}$.  Let $\mathcal{A}, \mathcal{X}$ be two sets. The indicator function yields $\mathbb{1}_{\mathcal{A}}(x)=1$ if $x\in {\mathcal{A}}$ and $\mathbb{1}_{\mathcal{A}}(x)=0$ otherwise. The set-wise difference of $\mathcal{A}, \mathcal{X}$ is $\mathcal{A}\setminus \mathcal{X} = \{x\in \mathcal{A}: x\notin  \mathcal{X}\}$. For a sequence $x_{0:N}$, we use $x_{0:N}\in \mathcal{A}$ to denote $x_k\in \mathcal{A},\: \forall \:k\in [N]$. 

\section{Preliminaries and Problem Formulation}
\label{sec_problem}
We define an MDP over a finite time horizon as a tuple $(\mathcal{X}, \mathcal{U}, T, \ell_{0:N})$, where the state space $\mathcal{X}$ and the input space $\mathcal{U}$ are Borel subsets of complete separable metric spaces equipped with the respective $\sigma$-algebras $\mathcal{B}(\mathcal{X})$ and $\mathcal{B}(\mathcal{U})$. Given a state $x_k\!\in\!\mathcal{X}$ and an input $u_k\!\in\!\mathcal{U}$, the stochastic kernel $T\!:\!\mathcal{B}(\mathcal{X})\!\times\!\mathcal{X}\!\times\!\mathcal{U}\!\rightarrow\![0,1]$ describes the stochastic state evolution, leading to $x_{k+1}\sim T(\cdot|x_k,u_k)$.
Additionally, $\ell_k\!:\!\mathcal{X}\!\times\mathcal{U}\!\rightarrow\!\mathbb{R}_{\geq 0}, k\in[N-1]$ and $\ell_N\!:\!\mathcal{X}\!\rightarrow\!\mathbb{R}_{\geq 0}$ denote measurable, non-negative functions called stage and terminal costs, respectively, which are incurred at every time step $k\in[N-1]$ and at terminal time $N\in\mathbb{N}$. 

 
A deterministic Markov policy is a sequence $\pi=(\mu_0,\mydots,\mu_{N-1})$ of measurable maps $\mu_k\!:\!\mathcal{X}\!\rightarrow\!\mathcal{U}, k\in [N-1]$. 
We denote the set of deterministic Markov policies by $\Pi$. In this work, we define mixed policies $\pi_{\text{mix}}$ as a deterministic Markov policy that is randomly chosen at the initial time step $k=0$ and then used during the entire control task. In other words, before applying any control inputs, a deterministic Markov policy is randomly selected from a set of deterministic Markov policies according to some probability law, and then this policy is used to compute control inputs for all time-steps. Note that this is in contrast to standard randomized policies, where an action is sampled at every time-step from a state or history dependent distribution. More formally, we endow the set $\Pi$ with the metric topology and define the corresponding Borel $\sigma$-algebra by $\mathcal{G}_{\Pi}$ \cite{aumann_1}. A mixed policy $\pi_{\text{mix}}$ can be seen as a random variable on the measurable space $(\Pi,\mathcal{G}_{\Pi})$. The set of all mixed policies is denoted as $\Pi_{\text{mix}}$. Note that $\Pi\subseteq\Pi_{\text{mix}}$.

For a given initial state $x_0\in\mathcal{X}$, policy $\pi\in\Pi$, and transition kernel $T$, a unique probability measure $\mathbb{P}_{x_0}^{\pi}$ for the state-input-trajectory is defined over $\mathcal{B}((\mathcal{X} \times \mathcal{U})^N\times \mathcal{X})$, which can be sampled recursively via $x_{k+1}\sim T(\cdot|x_k,\mu_k(x_k))$ \cite{abate2008probabilistic}. We define the associated expected cumulative cost as  
\begin{align}
    \begin{split}
            &\mathbb{E}_{x_0}^\pi 
            \left[ \ell_N(x_N) + \sum_{k=0}^{N-1} \ell_k(x_k,u_k) \right]=\hspace{-.4cm}\int\displaylimits_{(\mathcal{X} \times \mathcal{U})^N\times \mathcal{X}}\hspace{-.6cm}\biggl(\ell_N(x_N)\! \\ &\qquad+\!\sum_{k=0}^{N-1} \ell_k(x_k,u_k)\!\biggl)\mathbb{P}_{x_0}^{\pi}(dx_0,du_0,\dots,dx_N).
    \end{split}
\end{align}
Since the probability measure over trajectories $\mathbb{P}_{x_0}^{\pi}$ is unique, it follows that the expected cost is also unique.  

\textbf{A. Problem Formulation.}
In addition to minimizing cost we aim to synthesize controllers that meet additional specifications. Let $\mathcal{A},\mathcal{T}\in\mathcal{B}(\mathcal{X})$ define a safe set and target set, respectively, such that $\mathcal{A}\cap\mathcal{T}=\varnothing$. The set of trajectories satisfying specifications called invariance, reachability, and reach-avoidance, are defined respectively as
\begin{align}
        \mathcal{H}_{\text{I}}\!&=\! \{x_{0:N}\in\mathcal{X}^{N+1}\!:\!x_{0:N}\in\mathcal{A}\}, \nonumber \\
        \mathcal{H}_{\text{R}}\!&=\! \{x_{0:N}\in\mathcal{X}^{N+1}\!:\!\exists x_k\in x_{0:N}, x_k\in\mathcal{T}\}, \\
        \mathcal{H}_{\text{RA}}\!&=\! \{x_{0:N}\in\mathcal{X}^{N+1}\!:\!\exists x_k\in x_{0:N}, x_{0:{k-1}}\in\mathcal{A} \land x_k\in\mathcal{T}\}. \nonumber
\end{align}

Let $\mathcal{H}\in\{\mathcal{H}_{\text{I}},\mathcal{H}_{\text{R}},\mathcal{H}_{\text{RA}}\}$. Our aim is to minimize the cumulative cost, while guaranteeing a prescribed probability of satisfying the specification $\mathcal{H}$, leading to the joint chance constrained optimal control problem
\begin{equation}\label{eq_problemFormulation_CCOC}
\begin{aligned} 
    \begin{split}
    \inf_{\pi \in \Pi_{\text{mix}}} \quad  &\mathbb{E}_{x_0}^\pi\left[ \ell_N(x_N) + \sum_{k=0}^{N-1} \ell_k(x_k,u_k) \right]\\
    \text{s.t.} \quad  & \mathbb{P}_{x_0}^\pi(x_{0:N}\in \mathcal{H})\geq \alpha,
    \end{split}
\end{aligned}
\end{equation}
where $\alpha \in [0,1]$ is a user-specified success probability. 
As shown in \cite{schmid2023Jccoc}, the restriction in Problem \eqref{eq_problemFormulation_CCOC} to the class of mixed policies is without loss of optimality when compared to the class of stochastic causal policies.
For simplicity, we term the probability of satisfying the additional specification as safety probability, irrespective of the type of objective, and call a trajectory safe if $x_{0:N}\in\mathcal{H}$.

To ensure the existence of an optimal policy in Problem \eqref{eq_problemFormulation_CCOC}, we introduce the following standing assumption.
\begin{assumption}\label{assumption_of_attainability}
    The input set $\mathcal{U}$ is compact. For all $x\in\mathcal{X}$, $\mathcal{C}\in\mathcal{B}(\mathcal{X})$, the transition kernel $T(\mathcal{C}|x,u)$ and stage cost $\ell_k(x,u)$ are continuous with respect to $u$. A strictly feasible finite cost policy for Problem \eqref{eq_problemFormulation_CCOC} exists.
\end{assumption}
To simplify the notation we further assume that the stochastic kernel $T$ admits a density $t$. However, all of our results extend to kernels that admit no density. They also extend to finite state and action spaces. 



\textbf{B. Dynamic Programming.}
For a given policy $\pi\!\in\!\Pi$ and initial state $x_0\!\in\!\mathcal{X}$, the cost incurred over $N$ time steps is
$
    C_0^{\pi}(x_0)\!=\!\mathbb{E}_{x_{0}}^\pi\!\left[ \ell_N(x_N)\!+\!\sum_{k=0}^{N-1}\! \ell_k(x_k,\mu_k(x_k)) \right]
    \label{eq_min_DP_cost_problem}
$.
The cost can be recursively computed using the DP recursion \cite{bertsekas2012dynamic}
\begin{align}
        C_N^{\pi}(x_N)\!&=\!\ell_N(x_N),\label{cost_DP} \\
        C_k^{\pi}(x_k)\!&=\hspace{-0.2em}\ell_k(x_k,u_k)+\int_{\mathcal{X}}C_{k+1}^{\pi}(x_{k+1})T(dx_{k+1}|x_k,\!u_k),\nonumber
\end{align}
where $u_k\sim\mu_k(x_k)$.
The finite-horizon optimal control problem now aims to find the infimum expected cost over the policies $\pi\in\Pi$. Denoting $C_k^{\star}(x_k):=\inf_{\pi \in \Pi}C_k^{\pi}(x_k)$, we obtain \cite{bertsekas2012dynamic}
\begin{align}
        C_N^{\star}(x_N)\!&=\!\ell_N(x_N),\label{eq_max_cost_evaluation} \\
         C_k^{\star}(x_k)\!&=\hspace{-0.4em}\min_{u_k\in\mathcal{U}}\hspace{-0.2em}\ell(x_k,u_k)+\!\int_{\mathcal{X}} \! \!C_{k+1}^{\star}(x_{k+1})T(dx_{k+1}|x_k,u_k).\nonumber
\end{align}
The minimum is attained under Assumption \ref{assumption_of_attainability} \cite{ioffe2005generic}. The optimal deterministic Markov policy described by the DP recursion \eqref{eq_max_cost_evaluation} is also optimal within the class of mixed policies \cite[Theorem 3.2.1]{hernandez2012discrete}.

\textbf{C. Safety.}
Note that reachability specifications are a special case of reach-avoid specifications with the unsafe set of states being the empty set, i.e., $\mathcal{A}=\mathcal{X}\setminus\mathcal{T}$. Hence, we will drop the reachability case in the upcoming discussion. A similar reduction can be achieved for invariance specifications \cite{schmid2022probabilistic}. However, for clarity of presentation we will keep the invariance case. Success in safety tasks can be captured through products of set indicator functions as
\begin{align}
\begin{split}
    \prod_{k=0}^N\!\mathbb{1}_{\mathcal{A}}(x_k)&\!=\!\begin{cases}
        1 & \text{if $x_{0:N}\!\in\!\mathcal{H}_\text{I}$}, \\
        0 & \text{otherwise},
    \end{cases} \\
    1\!-\!\prod_{k=0}^N\!\left(\!1\!-\!\mathbb{1}_{\mathcal{T}}(x_k)\!\prod_{i=0}^{k-1}\!\mathbb{1}_{\mathcal{A}}(x_i)\! \right)\!&=\!\begin{cases}
        1 & \text{if $x_{0:N}\!\in\!\mathcal{H}_{\text{RA}}$}, \\
        0 & \text{otherwise}.
    \end{cases} 
\end{split}
\end{align}
The probability of success is then given as the expected value of the above terms; for instance, for invariance
$
    \mathbb{P}_{x_0}^{\pi}(x_{0:N}\in\mathcal{H}_{\text{I}}) = \mathbb{E}_{x_0}^{\pi}\left[\prod_{k=0}^N\mathbb{1}_{\mathcal{A}}(x_k)\right].
$
This quantity can be computed via a DP recursion with multiplicative cost, by recursively evaluating the product. We use $V_k^{\pi},V_k^{\star}:\mathcal{X}\rightarrow[0,1]$ as shorthand notation for $V_k^{\pi}(x_k)=\mathbb{P}_{x_0}^{\pi}(x_{k:N}\in \mathcal{A}|x_k, \pi)$ and $V_k^{\star}(x_k)=\sup_{\pi\in\Pi_{\text{mix}}} \mathbb{P}_{x_0}^{\pi}(x_{k:N}\in \mathcal{A}|x_k, \pi)$, respectively. Following  \cite{abate2008probabilistic}, for a given $\pi\in\Pi$ 
\begin{align}
    \begin{split}
        V_N^{\pi}(x_N) & =\mathbb{1}_{\mathcal{A}}(x_N), \\
        V_k^{\pi}(x_k) &= \mathbb{1}_{\mathcal{A}}(x_k)\hspace{-0.4em}\int_{\mathcal{X}}\hspace{-0.6em}V_{k+1}^\pi(x_{k+1})T(dx_{k+1}|x_{k}, u_k),
         \end{split}
    \label{eq_invariance_evaluation}
\end{align}
where $u_k\sim\mu_k(x_k)$, and
\begin{align}
    \begin{split}
        V_N^{\star}(x_N)\!&=\mathbb{1}_{\mathcal{A}}(x_N),\\
        V_k^{\star}(x_k)\!&=\!\max_{u_k\in\mathcal{U}}\hspace{-0.2em}\mathbb{1}_\mathcal{A}(x_k)\hspace{-0.4em}\int_{\mathcal{X}}\hspace{-0.6em}V_{k+1}^{\star}(x_{k+1})T(dx_{k+1}|x_{k}, u_k).
    \end{split}
    \label{eq_max_invariance_evaluation}
\end{align}
where the maximum is guaranteed to be attained by Assumption \ref{assumption_of_attainability} \cite{kariotoglou2017linear}. Similar recursions can be constructed to compute reachability and reach-avoid probabilities, see \cite{kariotoglou2017linear,schmid2022probabilistic}.

The existence of a feasible finite cost policy in Assumption \ref{assumption_of_attainability} can be verified by computing the safest policy, verifying that its safety is greater than $\alpha$ and its cost is finite via recursion \eqref{cost_DP}.

\section{Joint Chance Constrained DP}
\label{sec_jcc_dp}
We aim to solve Problem \eqref{eq_problemFormulation_CCOC} using its Lagrange-Dual
\begin{equation}\label{eq_dual_of_CCOC}
    \begin{split}
        \sup_{\lambda\in\mathbb{R}_{\geq 0}}\inf_{\pi \in \Pi_{\text{mix}}} \quad  &\mathbb{E}_{x_0}^\pi\left[ \ell_N(x_N) + \sum_{k=0}^{N-1} \ell_k(x_k,u_k) \right] \\ &+ \lambda(\alpha-\mathbb{P}_{x_0}^\pi(x_{0:N}\in \mathcal{H})).
    \end{split}
\end{equation}
However, the inner infimum is intractable due to the non-Markovian structure imposed by the constraint \cite{schmid2023Jccoc}. 

Recall that, by definition, $\mathcal{A}\cap\mathcal{T}=\varnothing$. For invariance specifications, we define a binary state $b_{\text{I},k}$, with $b_{\text{I},0} = \mathbb{1}_{\mathcal{A}}(x_0)$ and dynamics $b_{\text{I},k+1}=\mathbb{1}_{\mathcal{A}}(x_{k+1})b_{\text{I},k}$. Then, $\prod_{k=0}^N\mathbb{1}_{\mathcal{A}}(x_k) = b_{\text{I},N}$. Similarly, for reach-avoid specifications we define the variable $b_{\text{RA},k}\in\{0,1,2\}$, with $b_{\text{RA},0} = \mathbb{1}_{\mathcal{A}}(x_0) + 2\mathbb{1}_{\mathcal{T}}(x_0)$ and 
\begin{align}
    \hspace{-0.2cm}b_{\text{RA},k+1}\!=\!\begin{cases}
    2, & \text{if } b_{\text{RA},k}\!=\!2\lor(b_{\text{RA},k}\!=\!1 \land x_{k+1}\!\in\!\mathcal{T}),\\
    1, & \text{if } b_{\text{RA},k}\!=\!1 \land  x_{k+1}\!\in\!\mathcal{A}, \\
    0, & \text{otherwise.} 
    \end{cases}
\end{align}

\begin{figure}[!tbp]
    \centering
    \resizebox{\columnwidth}{!}{%
    \begin{tikzpicture}
        \node[state] (q1) {$b_{k}=1$};
        \node[state, right of=q1] (q2) {$b_k=0$};
        \draw (q1) edge[loop above] node{$x_{k+1}\in\mathcal{A}$} (q1)
        (q2) edge[loop above] node{} (q2)
        (q1) edge[bend left, above] node{$x_{k+1}\notin\mathcal{A}$} (q2);
        
        
        \node[state, right of=q2] (q6)  {$b_k=0$};
        \node[state, right of=q6](q5) {$b_k=1$};
        \node[state,right of=q5] (q7)  {$b_k=2$};
        \draw (q5) edge[loop above] node{$x_{k+1}\in\mathcal{A}$} (q5)
        (q6) edge[loop above] node{} (q6)
        (q7) edge[loop above] node{} (q7)
        (q5) edge[bend left=40, below] node{$x_{k+1}\notin\mathcal{A}\cup\mathcal{T}$} (q6)
        (q5) edge[bend right=40, below] node{$x_{k+1}\in\mathcal{T}$} (q7);
    \end{tikzpicture}
    }
    \caption{Augmented state dynamics for invariance (left) and reach-avoid objectives (right). Arrows indicate binary state transitions based on the annotated condition; unannotated arrows denote transitions under any condition.}
    \label{fig_binary_dynamics}
\end{figure}
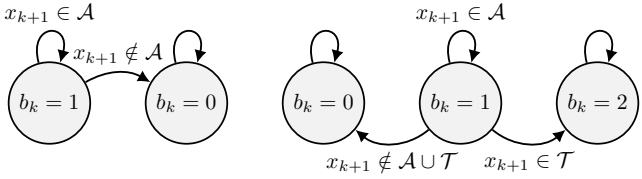

The dynamics are illustrated in Figure \ref{fig_binary_dynamics}.

\begin{lemma}
    The auxiliary state trajectories are such that
    \begin{align}
        b_{\text{I},N}=1 &\Leftrightarrow x_{0:N}\in\mathcal{H}_{\text{I}}, \ \ 
        b_{\text{RA},N}=2 \Leftrightarrow x_{0:N}\in\mathcal{H}_{\text{RA}}.
    \end{align}
\end{lemma}
\ifShowProofs
\begin{proof}
    To show the first equivalence, note that if there exists $k\in[N]$ such that $ x_k\notin\mathcal{A}$, then $b_{\text{I},k}=0$. By the stated dynamics, if $b_{\text{I},k}=0$, then $b_{\text{I},k+1}=0$, hence $b_{\text{I},N}=0$ by induction. Hence, if $b_{\text{I},N}=1$, then there does not exist $k\in[N]$ such that $x_k\notin\mathcal{A}$, hence $x_{0:N}\in\mathcal{H}_{\text{I}}$. Vice versa, if $x_{0:N}\in\mathcal{A}$, then $b_{\text{I},0}=1$ by $x_0\in\mathcal{A}$. Given $x_k\in\mathcal{A}$, $b_{\text{I},k}=1$, then $b_{\text{I},k+1}=1$; consequently $b_{\text{I},N}=1$ by induction. The proof for the reach-avoid case follows similar arguments and is omitted in the interest of space. 
\end{proof}
\fi

We augment the system dynamics by appending the auxiliary state to capture the success in the safety specifications over the trajectory. To simplify notation we drop the index I and RA from the variable $b_k$ in the sequel. Depending on the specification type, we define $
\Tilde{\mathcal{X}}= \mathcal{X}\times\{0,1\}$, or $
\Tilde{\mathcal{X}}= \mathcal{X}\times\{0,1,2\}$, and $\Tilde{x}_k = (x_k, b_k)\in\Tilde{\mathcal{X}}$. The combined dynamics of the MDP and auxiliary state are described by the stochastic kernels $\Tilde{T}_{\text{I}},\Tilde{T}_{\text{RA}}:\mathcal{B}(\Tilde{\mathcal{X}})\times\Tilde{\mathcal{X}}\times\mathcal{U}\rightarrow[0,1]$, with density
\begin{align}
    &\Tilde{t}_{\text{I}}(\Tilde{x}_{k+1}|\Tilde{x}_k,u_k)  
    \\&\!=\!\begin{cases}
        t(x_{k+1}|x_k,\!u_k)\!& \text{if} \ (b_{k}\!=\!0 \land b_{k+1}\!=\!0) 
        \\
        & \ \ \ \lor \  (b_{k}=1 \land b_{k+1}\!=\!\mathbb{1}_{\mathcal{A}}(x_{k+1})) 
        \\
        0 & \text{otherwise}, 
    \end{cases} \nonumber \\
    &\Tilde{t}_{\text{RA}}(\Tilde{x}_{k+1}|\Tilde{x}_k,u_k) \\&= \begin{cases}
        t(x_{k+1}|x_k,u_k) & \text{if} \ (b_{k}=2 \land b_{k+1}=2) \ \lor (b_k=1 \\ &  \ \ \ \land   \ b_{k+1}\!=\!\mathbb{1}_{\mathcal{A}}(x_{k+1})\!+\!2\mathbb{1}_{\mathcal{T}}(x_{k+1})) \\ & \ \ \ \lor \  (b_k=0 \land b_{k+1}=0) \\
        0 & \text{otherwise},
    \end{cases}\nonumber 
\end{align}
respectively, where $t$ is the density of $T$ and $\Tilde{x}_k = (x_k, b_k)$. Above dynamics implicitly comprise two steps. Firstly, the state component $x_k$ transitions to a state $x_{k+1}$ as described by $t$. Secondly, $b_{k+1}$ is deterministically computed based on $x_{k+1}$ and $b_k$. Combined, the density remains unchanged for combinations of $x_{k+1}$, $b_{k+1}$ and $b_k$ that are in accordance with the auxiliary state dynamics, and is set to zero otherwise. Reachability specifications can be treated as reach-avoid specifications. However, note that under reachability specifications a transition to a state with $b_k=0$ is impossible, allowing to reduce the state space to states with $b_k=1$ and $b_k=2$. 

We overload the cost functions $l_k(\Tilde{x}_k,u_k)=l_k(x_k,u_k),\:l_N(\Tilde{x}_N)=l_N(x_N), k\in[N-1]$. Lastly, we define success indicator functions
\begin{align}
    \delta_{\text{I}}(b_N)\!=\!\begin{cases}
    1 &\hspace{-0.3em}\text{if }b_N\!=\!1,\\ 
    0 &\hspace{-0.3em}\text{else, }
\end{cases} \  \delta_{\text{RA}}(b_N)\!=\!\begin{cases}
    1 &\hspace{-0.3em}\text{if }b_N\!=\!2,\\ 
    0 &\hspace{-0.3em}\text{else.}
\end{cases}
\end{align} 
For generality, we again drop the index from $\delta$ and $\Tilde{T}$ in the following and say that the trajectory is safe if $\delta(b_N)=1$. This allows to express the probability of safety by the expected value of the success indicator function, i.e., $\mathbb{P}_{x_0}^\pi(x_{0:N}\in \mathcal{H})=\mathbb{E}_{x_0}^\pi\left[\delta(b_N)\right]$; hence, we can rewrite Problem \eqref{eq_dual_of_CCOC} as
\begin{equation}\label{eq_dual_with_binary}
    \begin{split}
        \hspace{-0.58em}\sup_{\lambda\in\mathbb{R}_{\geq 0}}\hspace{-1.3em}\inf_{\ \ \  \pi \in \Pi_{\text{mix}}} \hspace{-0.85em}\mathbb{E}_{\Tilde{x}_0}^\pi\hspace{-0.45em}\left[\!\ell_N(\Tilde{x}_N)\!+\hspace{-0.5em}\sum_{k=0}^{N-1} \hspace{-0.2em}\ell_k(\Tilde{x}_k,u_k)\!+\!\lambda(\alpha\!-\!\delta(b_N))\!\right]
    \end{split}.
\end{equation}
Note that this problem formulation matches with \cite{schmid2023Jccoc}, replacing $b_N$ with $\delta(b_N)$. We now extend the results in \cite{schmid2023Jccoc}. 
\begin{theorem}
    The supremum and infimum in Problem \eqref{eq_dual_with_binary} are attained by some finite $\lambda \in \mathbb{R}_{\geq 0}$ and a mixed policy, respectively. The inner infimum in Problem \eqref{eq_dual_with_binary} is Lipschitz-continuous in $\lambda$. 
    \label{cor_max_min_attainable}
\end{theorem}
\ifUnshortenedProofs
\ifShowProofs
\ifArxivVersion
    \begin{proof}
        Given any fixed $\lambda\in\mathbb{R}_{\geq 0}$, the existence of an optimal policy for the inner infimum is guaranteed under Assumption \ref{assumption_of_attainability} \cite{schmid2023Jccoc}. Moreover, at least one such optimal policy is deterministic Markov \cite{schmid2023Jccoc}. This allows us to write $\min$ instead of $\inf$ in \eqref{eq_dual_with_binary} and restrict to the class $\Pi$ without the loss of optimality. 
        Given the existence of a strictly feasible finite cost policy (Assumption \ref{assumption_of_attainability}), an upper bound $\overline{\lambda}$ on $\lambda$ can be found \cite[Proposition 5.2]{schmid2023Jccoc}. Further, the function  $
            f(\lambda)\!=\!\min_{\pi \in \Pi}\mathbb{E}_{\Tilde{x}_0}^\pi\!\left[ \ell_N(\Tilde{x}_N)\!+\hspace{-0.4em}\sum_{k=0}^{N-1} \hspace{-0.2em}\ell_k(\Tilde{x}_k,u_k)\!+\!\lambda(\alpha\!-\!\delta(b_N))\right]$
        is Lipschitz-continuous in $\lambda$. To show this, consider any $\lambda_1,\lambda_2\in[0,\overline{\lambda}]$ and let $\pi_{\lambda_1},\pi_{\lambda_2}$ be respective minimizers in $f(\lambda_1)$ and $f(\lambda_2)$. Then, 
         $   f(\lambda_1)\!\leq\!\mathbb{E}_{\Tilde{x}_0}^{\pi_{{\lambda}_2}}\![ \ell_N(\Tilde{x}_N)\!+\hspace{-0.4em}\sum_{k=0}^{N-1} \hspace{-0.2em}\ell_k(\Tilde{x}_k,u_k)\!+\!\lambda_1(\alpha\!-\!\delta(b_N))]$.
        Further,
        \begin{align} 
            \vspace{-4em}
            f(\lambda_2) &+ |\lambda_1 - \lambda_2| \\&\geq f(\lambda_2) \\&=\mathbb{E}_{\Tilde{x}_0}^{\pi_{{\lambda}_2}}\!\Big[ \ell_N(\Tilde{x}_N)\!+\hspace{-0.4em}\sum_{k=0}^{N-1} \hspace{-0.2em}\ell_k(\Tilde{x}_k,u_k)\!+\!\lambda_1(\alpha-\delta(b_N)) \nonumber \\&\quad+ (\lambda_2 - \lambda_1)(\alpha-\delta(b_N))\Big] \\
            &= \mathbb{E}_{\Tilde{x}_0}^{\pi_{{\lambda}_2}}\!\Big[ \ell_N(\Tilde{x}_N)\!+\hspace{-0.4em}\sum_{k=0}^{N-1} \hspace{-0.2em}\ell_k(\Tilde{x}_k,u_k)\!+\!\lambda_1(\alpha-\delta(b_N))\Big]  \nonumber\\&\quad+ (\lambda_2 - \lambda_1)\mathbb{E}_{\Tilde{x}_0}^{\pi_{{\lambda}_2}}\!\left[\alpha-\delta(b_N)\right] \\
            &\geq f(\lambda_1) - |\lambda_1 - \lambda_2|,
        \end{align}
        where the first equality follows by $|\lambda_1 - \lambda_2|$ being non-negative, and the last inequality follows by the initially introduced inequalities and $-1\leq \alpha-\delta(b_N)\leq 1$. Reordering yields 
            $f(\lambda_1) - f(\lambda_2) \leq 2|\lambda_1 - \lambda_2|.$
        By symmetry, 
        $
            |f(\lambda_1) - f(\lambda_2)| \leq 2|\lambda_1 - \lambda_2|,
        $
        showing $f(\lambda)$ is Lipschitz-continuous in $\lambda$. Combined with the fact that $\lambda$ takes values in the compact set $[0,\overline{\lambda}]$, an optimal $\lambda$ is guaranteed to exist \cite{ioffe2005generic}. 
    \end{proof}
\else 
\begin{proof}
        Given any fixed $\lambda\in\mathbb{R}_{\geq 0}$, the existence of an optimal policy for the inner minimization is guaranteed under Assumption \ref{assumption_of_attainability} \cite{schmid2023Jccoc}. Moreover, at least one such optimal policy is deterministic Markov \cite{schmid2023Jccoc}. This allows us to write $\min$ instead of $\inf$ in \eqref{eq_dual_with_binary} and restrict minimization to the class $\Pi$. 
        
        Given the existence of a feasible finite cost policy (Assumption \ref{assumption_of_attainability}), bounds $[0,\overline{\lambda}]\ni\lambda$ can be found \cite{schmid2023Jccoc}. Further, the function 
        \begin{align*}
            f(\lambda)\!=\!\min_{\pi \in \Pi}\mathbb{E}_{\Tilde{x}_0}^\pi\!\left[ \ell_N(\Tilde{x}_N)\!+\hspace{-0.4em}\sum_{k=0}^{N-1} \hspace{-0.2em}\ell_k(\Tilde{x}_k,u_k)\!+\!\lambda(\alpha\!-\!\delta(b_N))\right]
        \end{align*}
        is Lipschitz-continuous in $\lambda$ (CITE ARXIV). 
        Combined, an optimal $\lambda$ is guaranteed to exist \cite{ioffe2005generic}. 
    \end{proof}
\fi 
\fi
\fi

Theorem \ref{cor_max_min_attainable} allows us to replace the supremum and infimum in Problem \eqref{eq_dual_with_binary} with a maximum and minimum. For a fixed $\lambda$, the inner minimization  forms an MDP that is optimally solved by a (not necessarily unique) deterministic Markov policy, obtained via the DP recursion \cite{schmid2023Jccoc}
\begin{align}
    J_N^{\lambda}(\Tilde{x}_N)&\!=\!\ell_N(\Tilde{x}_N) + \lambda(\alpha- \delta(b_N)),\label{eq_our_dp_recursion_innerdual} \\
    J_k^{\lambda}(\Tilde{x}_k)&\!=\!\min_{u_k\in\mathcal{U}}\!\ell_k(\Tilde{x}_k,u_k)\!+\hspace{-0.4em}\int_{\Tilde{\mathcal{X}}}\hspace{-0.6em}J_{k+1}^{\lambda}(\Tilde{x}_{k+1})\Tilde{T}(d\Tilde{x}_{k+1}|\Tilde{x}_k,u_k).
\nonumber
\end{align}
Then, any respective greedy policy on the value functions $J_{0:N}^{\lambda}$ is a minimizer in \eqref{eq_dual_with_binary} \cite{schmid2023Jccoc}. 

\begin{definition}
    Given $\lambda\in \mathbb{R}_{\geq 0}$, we denote as $\lambda$-optimal under $\Pi$ any policy $\pi_{\lambda}\in\Pi_{\lambda}$ with 
    \begin{align}
        \begin{split}
            \Pi_{\lambda}= \Big\{\pi'\in & \ \Pi: \pi'\in{\argmin}_{\pi\in\Pi}  \mathbb{E}_{\Tilde{x}_0}^\pi\big[\ell_N(\Tilde{x}_N) \\ &+ \sum_{k=0}^{N-1} \ell_k(\Tilde{x}_k,u_k)+\lambda (\alpha-\delta(b_N))\big]\Big\}.
        \end{split}
    \end{align}
\end{definition}


Let $\lambda^{\star}$ be an optimal argument for Problem \eqref{eq_dual_of_CCOC}, then $\Pi_{\lambda^{\star}}$ denotes the set of all $\lambda^{\star}$-optimal deterministic Markov policies for the inner minimization, which also yield optimal arguments in the class of mixed policies. Within $\Pi_{\lambda^{\star}}$, let $\pi_c$ be the policy with the lowest control cost and $\pi_v$ the policy with the highest probability on safety. Then, if Problem \eqref{eq_problemFormulation_CCOC} is feasible, based on \cite{schmid2023Jccoc} the following hold: First, $\pi_v$ is guaranteed to have a safety of at least $\alpha$ and second, an optimal mixed policy to Problem \eqref{eq_problemFormulation_CCOC} can be constructed from $\pi_c$ and $\pi_v$. In the next section, we will design an LP framework that follows this recipe. We first solve for an optimal dual multiplier $\lambda^{\star}$, then we compute the value functions $J_{0:N}^{\lambda^{\star}}$ which implicitly describe the set $\Pi_{\lambda^{\star}}$ and finally we solve for policies $\pi_c$ and $\pi_v$ to construct a mixed policy that solves Problem \eqref{eq_problemFormulation_CCOC} optimally.



\section{Joint Chance Constrained LP}
\label{sec_jcc_lp}
Consider the LP
\begin{equation}
    \begin{split}
        \sup_{\lambda\in\mathbb{R}_{\geq 0}, J_{0:N}\in \mathcal{M}_{\Tilde{\mathcal{X}}}}   \hspace{-1em} &J_0(\Tilde{x}_0)
        \\ \text{s.t.} \quad &  J_N(\Tilde{x}_N)\leq\ell_N(\Tilde{x}_N) + \lambda(\alpha- \delta(b_N)), \\
        &J_k(\Tilde{x}_k)\leq  \ell_k(\Tilde{x}_k,u_k)  \\ &\ \ + \int_{\Tilde{\mathcal{X}}} J_{k+1}(\Tilde{x}_{k+1})\Tilde{T}(d\Tilde{x}_{k+1}|\Tilde{x}_k,u_k),
    \end{split}
    \label{eq_lp_ccoc}
\end{equation}
where $\mathcal{M}_{\Tilde{\mathcal{X}}}$ denotes the set of measurable functions from $\Tilde{\mathcal{X}}$ to $\mathbb{R}$ and the constraints hold for all $x_k,x_N\in\Tilde{\mathcal{X}},u_k\in\mathcal{U},k\in[N-1]$. 
\begin{lemma}
    Let $\lambda\in\mathbb{R}_{\geq 0}$ be fixed in LP \eqref{eq_lp_ccoc}. Any feasible $J_{0:N}^{\text{\eqref{eq_lp_ccoc}}}$ to LP \eqref{eq_lp_ccoc} yields $J_k^{\text{\eqref{eq_lp_ccoc}}}(\Tilde{x}_k)\leq J_{k}^{{\lambda}}(\Tilde{x}_k)$ for all $k\in[N]$, $\Tilde{x}_k\in\Tilde{\mathcal{X}}$, where $J_{0:N}^{{\lambda}}$ is defined by recursion \eqref{eq_our_dp_recursion_innerdual}. 
    \label{lem_lp_lower_bounds_dp}
\end{lemma}
\ifUnshortenedProofs
\ifShowProofs
\begin{proof}
    Since $J_N^{\text{\eqref{eq_lp_ccoc}}}$ is feasible for LP \eqref{eq_lp_ccoc} and by recursion \eqref{eq_our_dp_recursion_innerdual} we have for every $\Tilde{x}_N\in\Tilde{\mathcal{X}}$
    \begin{align}
        J_N^{\text{\eqref{eq_lp_ccoc}}}(\Tilde{x}_N) \leq \ell_N(\Tilde{x}_N) + \lambda(\alpha-\delta(b_N)) = J_N^{\lambda}(\Tilde{x}_N).
    \end{align}
    Further, suppose that there exists $k\in[N-1],\Tilde{x}_k\in\Tilde{\mathcal{X}}$ such that $J_k^{\text{\eqref{eq_lp_ccoc}}}(\Tilde{x}_k)> J_k^{\lambda}(\Tilde{x}_k)$. Then, for some input $u_k\in\mathcal{U}$, 
    \begin{align}
        \hspace{-0.3cm}J_k^{\lambda}(\Tilde{x}_k)&\!=\!\ell_k(\Tilde{x}_k,u_k)\!+\hspace{-0.15cm}\int_{\Tilde{\mathcal{X}}}\hspace{-0.15cm}J_{k+1}^{\lambda}(\Tilde{x}_{k+1})\Tilde{T}(d\Tilde{x}_{k+1}|\Tilde{x}_k,u_k)\\ &\!<\!J_k^{\text{\eqref{eq_lp_ccoc}}}(\Tilde{x}_k)
        \\&\!\leq\!\ell_k(\Tilde{x}_k,u_k)\!+\hspace{-0.15cm}\int_{\Tilde{\mathcal{X}}}\hspace{-0.15cm}J^{\text{\eqref{eq_lp_ccoc}}}_{k+1}(\Tilde{x}_{k+1})\Tilde{T}(d\Tilde{x}_{k+1}|\Tilde{x}_k,u_k).
    \end{align}  Since $\Tilde{T}(\cdot|\cdot,\cdot)$ is non-negative, there must exist $\Tilde{x}_{k+1}\in\Tilde{\mathcal{X}}$ where $J_{k+1}^{\lambda}(\Tilde{x}_{k+1})< J_{k+1}^{\text{\eqref{eq_lp_ccoc}}}(\Tilde{x}_{k+1})$. Proceeding inductively up to time step $N$ we obtain that there must exist $\Tilde{x}_N\in\Tilde{\mathcal{X}}$ such that
    $J_N^{\lambda}(\Tilde{x}_N) = \ell_N(\Tilde{x}_N) + \lambda(\alpha-\delta(b_N)) < J_N^{\text{\eqref{eq_lp_ccoc}}}(\Tilde{x}_N)\leq \ell_N(\Tilde{x}_N) + \lambda(\alpha-\delta(b_N)),$
    which is a contradiction. Hence, $J_k^{\text{\eqref{eq_lp_ccoc}}}(\Tilde{x}_k)\leq J_k^{\lambda}(\Tilde{x}_k)$ for all $k\in[N],\Tilde{x}_k\in\Tilde{\mathcal{X}}$.
\end{proof}
\fi
\fi
\begin{theorem}
    The supremum in LP \eqref{eq_lp_ccoc} is attained by some $\lambda\in\mathbb{R}_{\geq 0}$, $J_{0:N}^{\text{\eqref{eq_lp_ccoc}}}\in \mathcal{M}_{\Tilde{\mathcal{X}}}$. An optimal argument $\lambda$ of the LP \eqref{eq_lp_ccoc} is also optimal for Problem \eqref{eq_dual_of_CCOC} and vice versa.
    \label{thm_optimal_lambda}
\end{theorem}
\ifShowProofs
\begin{proof}
    Fix any $\lambda\in\mathbb{R}_{\geq 0}$. Then, by Lemma \ref{lem_lp_lower_bounds_dp} any feasible solution of the LP yields a lower bound to the DP solution $J_{0:N}^{{\lambda}}(\cdot)$. Further, $J_{0:N}^{{\lambda}}(\cdot)$ are measurable functions by \cite{schmid2023Jccoc} and satisfy the constraints in LP \eqref{eq_lp_ccoc}. Hence, under the fixed $\lambda$, they are an optimal solution of the LP. 

    Let $\Lambda^{\star}$ denote the set of optimal values for $\lambda$ in Problem \eqref{eq_dual_of_CCOC}, i.e., the maximum value of $J_{0}^{\lambda}(\Tilde{x}_0)$ is achieved if and only if $\lambda\in\Lambda^{\star}$. 
    Fix any $\lambda'\notin\Lambda^{\star}$ in LP \eqref{eq_lp_ccoc} and let $\lambda^{\star}\in\Lambda^{\star}$. Then, $J_{0}^{\text{\eqref{eq_lp_ccoc}}}(\Tilde{x}_0)\leq J_{0}^{\lambda'}(\Tilde{x}_0)$ by Lemma \ref{lem_lp_lower_bounds_dp} and $J_{0}^{\lambda'}(\Tilde{x}_0)< J_{0}^{\lambda^{\star}}(\Tilde{x}_0)$. Combined, $J_{0}^{\text{\eqref{eq_lp_ccoc}}}(\Tilde{x}_0)< J_{0}^{\lambda^{\star}}(\Tilde{x}_0)$ and together with $(\lambda^{\star},J_{0:N}^{\lambda^{\star}})$ being feasible in LP \eqref{eq_dual_of_CCOC} $J_{0:N}^{\text{\eqref{eq_lp_ccoc}}}$ is a suboptimal solution to LP \eqref{eq_dual_of_CCOC}. Hence, any optimal $\lambda$ in \eqref{eq_lp_ccoc} must be contained in $\Lambda^{\star}$. Vice versa, any $\lambda^{\star}\in\Lambda^{\star}$ is optimal to LP \eqref{eq_lp_ccoc} since $(\lambda^{\star},J_{0:N}^{\lambda^{\star}})$ is a feasible and optimal solution that attains the greatest possible value of $J_{0}^{\lambda^{\star}}(\Tilde{x}_0)$.
\end{proof}
\fi

While the above Theorem guarantees that the optimal arguments for $\lambda$ in LP \eqref{eq_lp_ccoc} and Problem \eqref{eq_dual_of_CCOC} are equivalent, this is not necessarily true for the value functions $J_{0:N}^{\text{\eqref{eq_lp_ccoc}}}$ and $J_{0:N}^{\lambda}$. Hence, greedy policies for $J_{0:N}^{\text{\eqref{eq_lp_ccoc}}}$ advertise potentially suboptimal inputs, making $J_{0:N}^{\text{\eqref{eq_lp_ccoc}}}$ useless for constructing an optimal solution for Problem \eqref{eq_problemFormulation_CCOC} \cite{schmid2022probabilistic}. We address this issue by solving a second LP, keeping $\lambda$ fixed to its optimal value $\lambda^{\star}$ from LP \eqref{eq_lp_ccoc} and including the value functions $J_{0:N}$ in the objective via a non-negative measure $\nu$ that assigns positive
mass to all open subsets of $\Tilde{\mathcal{X}}$:
\begin{equation}
    \begin{split}
        \max_{J_{0:N}\in \mathcal{M}_{\Tilde{\mathcal{X}}}} \quad  &J_0(\Tilde{x}_0)+\sum_{k=1}^N\int_{\Tilde{\mathcal{X}}}J_k(\Tilde{x}_k)\nu(dx_k)
        \\ \text{s.t.} \quad &J_N(\Tilde{x}_N)\leq\ell_N(\Tilde{x}_N) + \lambda^{\star}(\alpha- \delta(b_N)), \\
        &J_k(\Tilde{x}_k)\leq  \ell_k(\Tilde{x}_k,u_k) + \\ & \ \quad \qquad \int_{\Tilde{\mathcal{X}}} J_{k+1}(\Tilde{x}_{k+1})\Tilde{T}(d\Tilde{x}_{k+1}|\Tilde{x}_k,u_k),
    \end{split}\label{eq_lp_ccoc_value_functions}
\end{equation}
where the constraints hold for all $\Tilde{x}_k,\Tilde{x}_N\in\Tilde{\mathcal{X}},u_k\in\mathcal{U},k\in[N-1]$. We denote the optimal arguments as $J_{0:N}^{\text{\eqref{eq_lp_ccoc_value_functions}}}$. The existence of the maximum follows as in Theorem \ref{thm_optimal_lambda}. 
\begin{theorem}
\label{thm_lp_solution_optimal}
    Let $\lambda^{\star}$ denote the optimal argument of $\lambda$ in LP \eqref{eq_lp_ccoc}. The value functions $J_{0:N}^{\text{\eqref{eq_lp_ccoc_value_functions}}}$ obtained as the optimal argument of LP \eqref{eq_lp_ccoc_value_functions} are upper bounded by $J_{0:N}^{\lambda^{\star}}$ obtained by the DP recursion \eqref{eq_dual_of_CCOC} and are $\nu$-almost-everywhere equivalent. The values $J_{0}^{\text{\eqref{eq_lp_ccoc_value_functions}}}(\Tilde{x}_0) = J_{0}^{\lambda^{\star}}(\Tilde{x}_0)$ are equivalent.
\end{theorem}
\ifUnshortenedProofs
\ifShowProofs
\begin{proof}
    The fact $J_{0:N}^{\text{\eqref{eq_lp_ccoc_value_functions}}}\leq J_{0:N}^{\lambda^{\star}}$ can be established as in Lemma \ref{lem_lp_lower_bounds_dp} by $J_{0:N}^{\lambda^{\star}}$ being feasible for LP \eqref{eq_lp_ccoc_value_functions}. If there is a non-zero $\nu$-measure set of states where the LP solution is strictly suboptimal or if $J_0^{\text{\eqref{eq_lp_ccoc_value_functions}}}(\Tilde{x}_0)$ is strictly suboptimal, then 
    $        J_0^{\text{\eqref{eq_lp_ccoc_value_functions}}}(\Tilde{x}_0)\!+\!\sum_{k=1}^N\!\int_{\Tilde{\mathcal{X}}}\!J_k^{\text{\eqref{eq_lp_ccoc_value_functions}}}(\Tilde{x}_k)\nu(d\Tilde{x}_k)\!<\!J_0^{\lambda^{\star}}(\Tilde{x}_0)\!+\!\sum_{k=1}^N\!\int_{\Tilde{\mathcal{X}}}\!J_k^{\lambda^{\star}}(\Tilde{x}_k)\nu(d\Tilde{x}_k)
        $
    and, by feasibility of $J_{0:N}^{\lambda^{\star}}$, $J_{0:N}^{\text{\eqref{eq_lp_ccoc_value_functions}}}$ is not an optimal solution of the LP. Hence, $J_{0:N}^{\text{\eqref{eq_lp_ccoc_value_functions}}}$ can not be smaller than  $J_{0:N}^{\lambda^{\star}}$ for a non-zero $\nu$-measure set of states and not larger by the prior bound, leading to the claim.
\end{proof}
\fi
\fi


\begin{assumption}
    For every $\Tilde{x}\in\Tilde{\mathcal{X}},u\in\mathcal{U}$, the kernel $\Tilde{T}(\cdot|\Tilde{x},u)$ is absolutely continuous with respect to $\nu$. 
    \label{ass_T_abs_cont}
\end{assumption}
\begin{proposition}
    Under Assumption \ref{ass_T_abs_cont}, the set of greedy policies with respect to $J_{0:N}^{\text{\eqref{eq_lp_ccoc_value_functions}}}(\Tilde{x}_0)$ is equivalent to the set of greedy policies with respect to $J_{0:N}^{\lambda^{\star}}(\Tilde{x}_0)$.
    \label{prop_value_functions_of_dp_and_lp_equal}
\end{proposition}
\ifUnshortenedProofs
\ifShowProofs
\begin{proof}
    For $k\in[N]$, let $\mathcal{C}_k=\{\Tilde{x}\in\Tilde{\mathcal{X}}|J_{k}^{\text{\eqref{eq_lp_ccoc_value_functions}}}(\Tilde{x}_k)< J_{k}^{\lambda^{\star}}(\Tilde{x}_k)\}$, which has zero $\nu$-measure by \cref{thm_lp_solution_optimal}. Then, for $k\in[N-1]$ the greedy policy picks 
    \begin{align}
        &\argmin_{u_k\in\mathcal{U}}\int_{\Tilde{X}}J_{k+1}^{\text{\eqref{eq_lp_ccoc_value_functions}}}(\Tilde{x}_{k+1})\Tilde{T}(d\Tilde{x}_{k+1}|\Tilde{x}_k,u_k)
        \\=&\argmin_{u_k\in\mathcal{U}}\int_{\Tilde{X}\setminus\mathcal{C}_{k+1}}\hspace{-2em}J_{k+1}^{\text{\eqref{eq_lp_ccoc_value_functions}}}(\Tilde{x}_{k+1})\Tilde{T}(d\Tilde{x}_{k+1}|\Tilde{x}_k,u_k) 
    \\=&\argmin_{u_k\in\mathcal{U}}\int_{\Tilde{X}\setminus\mathcal{C}_{k+1}}\hspace{-2em}J_{k+1}^{\lambda^{\star}}(\Tilde{x}_{k+1})\Tilde{T}(d\Tilde{x}_{k+1}|\Tilde{x}_k,u_k)
        \\=&\argmin_{u_k\in\mathcal{U}}\int_{\Tilde{X}}J_{k+1}^{\lambda^{\star}}(\Tilde{x}_{k+1})\Tilde{T}(d\Tilde{x}_{k+1}|\Tilde{x}_k,u_k),
    \end{align} since $\mathcal{C}_{k+1}$ is associated zero $\nu$-, hence, zero $\Tilde{T}$-measure.
    %
\end{proof}
\fi
\fi

In summary, solving the LP \eqref{eq_lp_ccoc} returns the optimal dual multiplier $\lambda^{\star}$. Subsequently solving the LP \eqref{eq_lp_ccoc_value_functions} under $\lambda^{\star}$ returns the value functions $J_{0:N}^{\text{\eqref{eq_lp_ccoc_value_functions}}}$, whose greedy policies correspond to the set of $\lambda^{\star}$-optimal policies. Following \cite{schmid2023Jccoc}, the optimal mixed policy for Problem \eqref{eq_problemFormulation_CCOC} can now be constructed by interpolating the cheapest and safest $\lambda^{\star}$-optimal deterministic Markov policies, denoted $\pi_{c,\lambda^{\star}}$, $\pi_{v,\lambda^{\star}}$, respectively. These can be computed using the recursions \eqref{eq_max_cost_evaluation} and \eqref{eq_max_invariance_evaluation}, restricting to inputs that are also optimal under $J_{0:N}^{\text{\eqref{eq_lp_ccoc_value_functions}}}$, e.g., the safest $\lambda^{\star}$-optimal policy can be computed with the recursion 
\begin{align}
        V_{v,N}^{\star}(\Tilde{x}_N) \!&=\hspace{-0.1em}\mathbb{1}_{\mathcal{A}}(\Tilde{x}_N), \label{eq_safety_evaluation_of_lambda_optimal} \\
        V_{v,k}^{\star}(\Tilde{x}_k)\!&=\hspace{-0.4em}\sup_{u_k\in\mathcal{U}} \mathbb{1}_\mathcal{A}(\Tilde{x}_k)\int_{\Tilde{\mathcal{X}}}\hspace{-0.4em}V_{v,k+1}^{\star}(\Tilde{x}_{k+1})\Tilde{T}(d\Tilde{x}_{k+1}|\Tilde{x}_{k}, u_k), \nonumber\\
        &\ \hspace{0.8em}\text{s.t.} \int_{\Tilde{X}}\hspace{-0.4em}J_{k+1}^{\text{\eqref{eq_lp_ccoc_value_functions}}}(\Tilde{x}_{k+1})\Tilde{T}(d\Tilde{x}_{k+1}|\Tilde{x}_k,u_k) \nonumber \\&\ \hspace{0.0em}= \argmin_{u'\in\mathcal{U}}\int_{\Tilde{X}}\hspace{-0.4em}J_{k+1}^{\text{\eqref{eq_lp_ccoc_value_functions}}}(\Tilde{x}_{k+1})\Tilde{T}(d\Tilde{x}_{k+1}|\Tilde{x}_k,u'),\nonumber
\end{align}
for $k\in[N-1]$. The constraint is continuous in $u_k$ since $\Tilde{T}$ is continuous in $u_k$ given a fixed $\Tilde{x}_k$. Hence, the level set described by the constraint is closed. Intersected with the constraints described by the compact set $\mathcal{U}$, the constraints are compact. Since the objective function is continuous in $u_k$ given $\Tilde{x}_k$, the supremum is always attained by some input $u_k$ \cite{ioffe2005generic}. A similar recursion can be formulated based on \eqref{eq_max_cost_evaluation} to compute the cheapest $\lambda^{\star}$-optimal policy $\pi_{c,\lambda^{\star}}$.

\begin{proposition}
    Under Assumption \ref{assumption_of_attainability} and \ref{ass_T_abs_cont}, let $\lambda^{\star}$ be an optimizer of \eqref{eq_lp_ccoc} and $J_{0:N}^{\eqref{eq_lp_ccoc_value_functions}}$ a maximizer of \eqref{eq_lp_ccoc_value_functions}. Let $\pi_{c,\lambda^{\star}}$ and $\pi_{v,\lambda^{\star}}$ be the cheapest and safest greedy policies on $J_{0:N}^{\eqref{eq_lp_ccoc_value_functions}}$ with associated safety $v_{c}$ and $v_{v}$, respectively. Then, 
    \begin{enumerate}
        \item[(i)] the policy $\pi_{v,\lambda^{\star}}$ is feasible in \eqref{eq_problemFormulation_CCOC};
        \item[(ii)] if $v_{c}$ is greater or equal to $\alpha$, then $\pi_{c,\lambda^{\star}}$ is an optimal solution to Problem \eqref{eq_problemFormulation_CCOC}. Otherwise, an optimal policy to Problem \eqref{eq_problemFormulation_CCOC} is obtained by sampling
        \vspace{-1mm}
    \begin{align}
    \pi_{\text{mix}}=\begin{cases}
    \pi_{v,\lambda^{\star}} & \text{with probability } \frac{\alpha - v_{c}}{v_{v} - v_{c}}, \\
    \pi_{c,\lambda^{\star}} & \, \text{otherwise.}
    \end{cases}
    \label{eq_stochastic_policies_ratios}
    \end{align}
    \end{enumerate}
\end{proposition}
\ifShowProofs
\vspace{1mm}
\begin{proof}
    By definition, any $\lambda^{\star}$-optimal policy solves the dual problem \eqref{eq_dual_with_binary} and any $\lambda^{\star}$-optimal mixed policy is constructed purely from $\lambda^{\star}$-optimal policies \cite[Corollary 4.7]{schmid2023Jccoc}. There exists a $\lambda^{\star}$-optimal mixed policy, which solves Problem \eqref{eq_problemFormulation_CCOC}; its existence is guaranteed under Assumption \ref{assumption_of_attainability} and its safety and cost is obtained by mixing at most two $\lambda^{\star}$-optimal deterministic Markov policies \cite[Theorem 4.10, Corollary 4.5]{schmid2023Jccoc}. Consequently, there must exist at least one $\lambda^{\star}$-optimal deterministic Markov policy with a safety of $\alpha$, in particular $\pi_{v,\lambda^{\star}}$. 
    
    Furthermore, the optimal mixed policy to Problem \eqref{eq_problemFormulation_CCOC} is that, which has the lowest cost, but a safety of at least $\alpha$. By the monotonicity property \cite[Lemma 4.8]{schmid2023Jccoc}, if $\pi_{c,\lambda^{\star}}$ has a safety of at least $\alpha$, it will be optimal. Otherwise, mixing any two $\lambda^{\star}$-optimal policies, e.g., $\pi_{c,\lambda^{\star}}$, $\pi_{v,\lambda^{\star}}$ via $\eqref{eq_stochastic_policies_ratios}$, to yield a mixed policy with safety $\alpha$ yields an optimal policy for Problem \eqref{eq_problemFormulation_CCOC}. No other $\lambda^{\star}$-optimal mixed policy can have a lower cost with at least the same safety, otherwise this policy would also be superior in the dual objective \eqref{eq_dual_with_binary}. This is a contradiction by $\lambda^{\star}$-optimality of the mixed policy $\eqref{eq_stochastic_policies_ratios}$. 
\end{proof}
\fi

We summarize the algorithmic steps in \cref{alg_final}.

\begin{algorithm}[!tbp]
\caption{Joint Chance Constr. Optimal Control}
$\lambda^{\star}\leftarrow$ LP \eqref{eq_lp_ccoc}\;
$J_{0:N}^{\eqref{eq_lp_ccoc_value_functions}} \leftarrow$ LP \eqref{eq_lp_ccoc_value_functions}\;
$(\pi_{c,\lambda^{\star}}) \leftarrow$ Recursion \eqref{eq_max_cost_evaluation} s.t. greedy on $J_{0:N}^{\text{\eqref{eq_lp_ccoc_value_functions}}}$\;
$v_{c} \leftarrow$ Recursion \eqref{eq_invariance_evaluation} s.t. $\pi_{c,\lambda^{\star}}$\;
$(\pi_{v,\lambda^{\star}},v_{v}) \leftarrow$ Recursion \eqref{eq_max_invariance_evaluation} s.t. greedy on $J_{0:N}^{\text{\eqref{eq_lp_ccoc_value_functions}}}$\;
$\pi_{\text{mix}}\leftarrow$ Eq. \eqref{eq_stochastic_policies_ratios}\;
return $\pi_{\text{mix}}$\;
\label{alg_final}
\end{algorithm}

\section{Numerical Examples}
\label{sec_examples}
\begin{figure}
    \centering
        \begin{subfigure}[t]{0.33\columnwidth}
            \centering
            \includegraphics[width=\columnwidth]{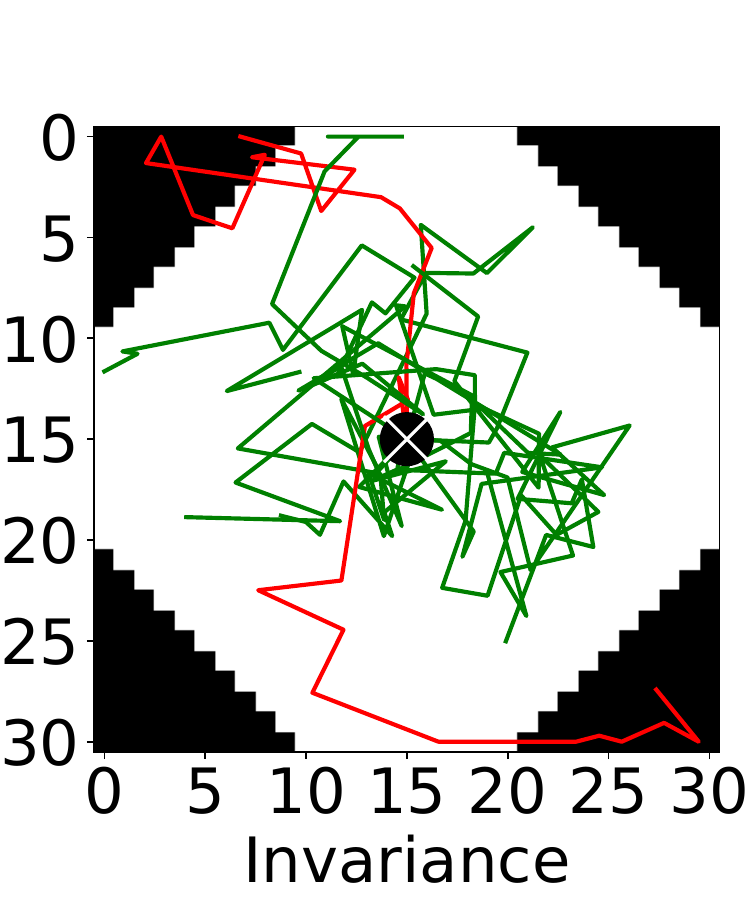}
        \end{subfigure}%
        \hfill
        \begin{subfigure}[t]{0.33\columnwidth}
            \centering
            \includegraphics[width=\columnwidth]{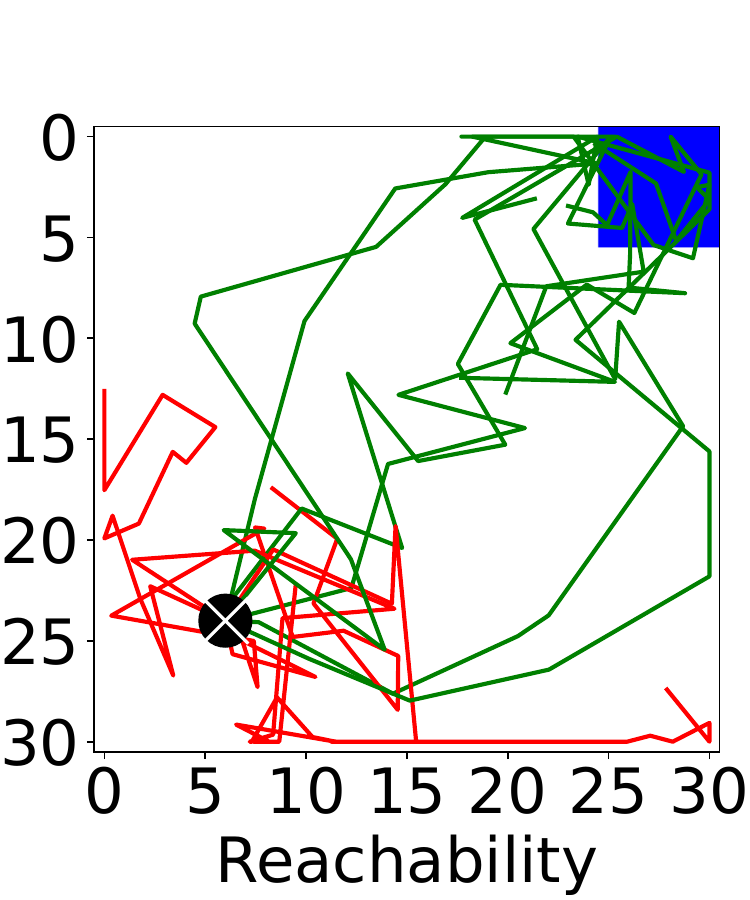}
        \end{subfigure}%
        \hfill
        \begin{subfigure}[t]{0.33\columnwidth}
            \centering
            \includegraphics[width=\columnwidth]{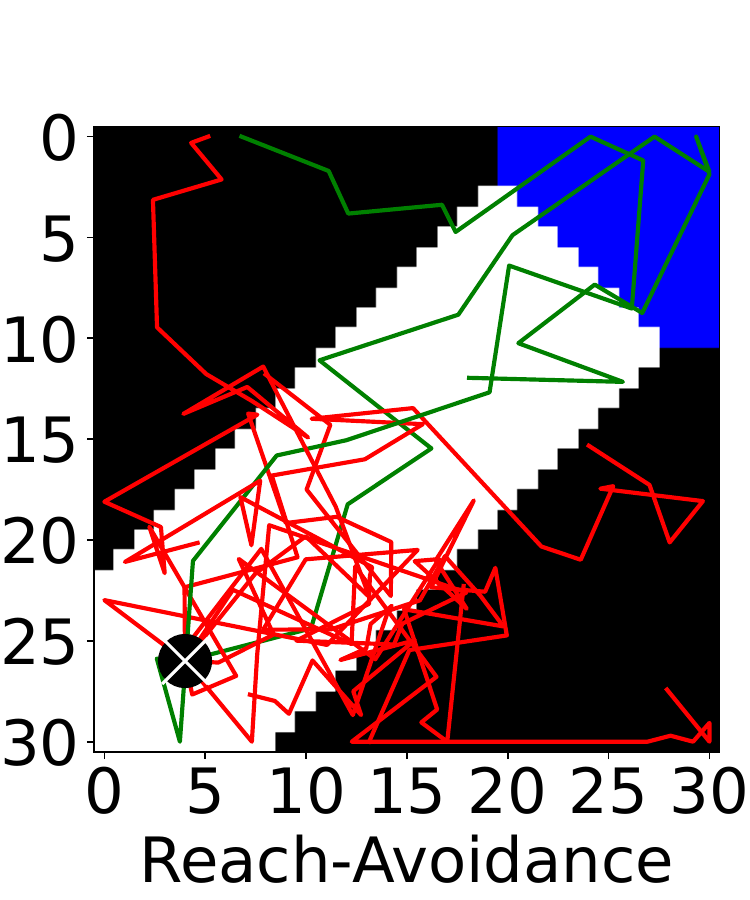}
        \end{subfigure}%
    \caption{The plots show ten trajectories generated by the optimal mixed policies for the invariance, reachability and reach-avoid problem in Section \ref{sec_examples}. The black regions denote unsafe and the blue regions target states. Trajectories that achieve the objective are painted green, and red otherwise. The crossed circle marks the initial state $x_0$.} 
    \label{fig_simulation_plots}
\end{figure}

To illustrate the theoretical results we consider a unicycle model with state space $\mathcal{X}=[0,10]\times[0,10]$, input space $\mathcal{U}=[0,3]\times[0,2\pi]$, dynamics $x_{k+1} = x_k + u_{1,k}\begin{bmatrix}\cos(v_k) & \sin(v_k)\end{bmatrix}^{\top} + w_{1,k}$, where $v_k = u_{2,k} + w_{2,k}$, with $u_k=\begin{bmatrix}
    u_{1,k} & u_{2,k}
\end{bmatrix}^{\top}$ representing the input and $w_{1,k}=\mathcal{N}([0 \ 0]^{\top},\text{diag}(1,1)),\:w_{2,k}=\mathcal{N}(0,0.2)$ the disturbances. The cost is defined by $\ell_k(x_k,u_k)=u_{1,k}, \ell_N(x_N)=0$ for all $x_k,x_N\in\mathcal{X},u_k\in\mathcal{U},k\in[N-1]$. Since the state and actions space are continuous there exist infinitely many states and inputs, resulting in LP \eqref{eq_lp_ccoc} and \eqref{eq_lp_ccoc_value_functions} to be infinite dimensional. As a counter-measure, we approximate the system dynamics via gridding, leading to finitely many states and actions \cite{schmid2022probabilistic}, but also other basis function approximations could have been utilized \cite{de2003linear}. 
The state and action space are discretized into $11\times 11$ states and $3\times 4$ actions respectively and the transition probabilities from each state under every input are estimated using $400$ Monte-Carlo trials. We consider an invariance, reachability and reach-avoid task over $N=15$ time steps. An optimal joint chance constrained control policy is computed for all scenarios using Algorithm \ref{alg_final} under the $\alpha$ reported in Table \ref{table_sim_results_numbers}. Ten simulated trajectories under the respectively optimal mixed policies applied to the continuous dynamics are presented in Figure \ref{fig_simulation_plots}. Table \ref{table_sim_results_numbers} further lists the cost and safety of the maximum safe, minimum cost, maximum safe $\lambda$-optimal, minimum cost $\lambda$-optimal and respective mixed policy when evaluated on the discretized dynamics and when applied to the continuous dynamics. To retrieve an estimate for the latter a Monte-Carlo simulation is performed using $10.000$ trials. Running algorithm \ref{alg_final} took $391$ seconds for the invariance and reachability problem and $579$ seconds for the reach-avoid problem using a Ryzen 5 3600 processor at 3.6 GHz and 16 GB DDR4 RAM at 3.6 GHz, CL16. The code used to generate the results is available at \url{https://github.com/NiklasSchmidResearch/JCC_opt_control_LP}.

The optimal mixed policy meets the desired safety specification precisely under the discretized dynamics in every scenario. When applied to the continuous dynamics, the safety probability is $1-5\%$ lower in all settings, which indicates errors introduced by the gridding approximation. However, compared to the safest policy, the optimal mixed policy yields a significantly reduced cost when evaluated on the discretized as well as when applied to the continuous dynamics. For instance, the safest policy for the invariance problem yields a cost of $11.87$ under the discretized ($12.03$ under the continuous) dynamics at a safety of $94.70\% (92.72 \%)$, while the optimal mixed policy yields a cost of $3.63 (3.94)$ at a safety of $90.00 \% (86.35\%)$. Significant cost reductions can also be observed for the reachability and reach-avoid example (see Table \ref{table_sim_results_numbers}). 

\begin{table}[]
    \begin{tabular}{llll}
    Specification               & Invariance     & Reachability   & Reach-Avoid.   \\ \hline
    $\alpha$                    & 0.90           & 0.60           & 0.25           \\
    $\pi_v$ (discr.)            & (11.87, 94.70) & (24.87, 88.24) & (12.20, 26.82) \\
    $\pi_c$ (discr.)            & (0.00, 58.67)  & (0.00, 0.32)   & (0.00, 0.14)   \\
    $\pi_{v,\lambda}$ (discr.)  & (3.63, 90.00)  & (13.57, 64.04) & (9.93, 25.66) \\
    $\pi_{c,\lambda}$ (discr.)  & (3.63, 90.00)  & (10.45, 51.72) & (9.30, 24.68) \\
    $\pi_{\text{mix}}$ (discr.) & (3.63, 90.00)  & (12.55, 60.00) & (9.51, 25.00) \\
    $\pi_{\text{mix}}$ (cont.)   & (3.95, 86.35)  & (12.11, 55.81) & (9.39, 23.40)
    \end{tabular}
       \caption{The table lists the safety specification $\alpha$, as well as the cost and safety (in $\%$) of five different policies when applied to the discretized (discr.) and continuous dynamics (cont.) for the invariance, reachability and reach-avoid problems in Section \ref{sec_examples}. The value of $\lambda$ is the result of the LP \eqref{eq_lp_ccoc}. The policies are $\pi_v$: The safest policy (recursion \eqref{eq_max_invariance_evaluation}), $\pi_c$: The lowest control cost policy (recursion \eqref{eq_max_cost_evaluation}), $\pi_{v,\lambda}$: The $\lambda$-optimal policy with the highest safety (recursion \eqref{eq_safety_evaluation_of_lambda_optimal}), $\pi_{c,\lambda}$: The $\lambda$-optimal policy with the lowest control cost, and $\pi_{\text{mix}}$: the optimal mixing of $\pi_{v,\lambda}$ and $\pi_{c,\lambda}$.}
  \label{table_sim_results_numbers}
    \end{table}

{
\section{Conclusion}
\label{sec:conclusion}

We presented an LP based framework to solve joint chance constrained optimal control problems with invariance, reachability or reach-avoid specifications.  
Multiple extensions are possible. Our MDP-based approach suffers from the curse of dimensionality. Basis function approximations, applicable to our LPs, might improve the scalability \cite{de2003linear}. More general specifications can be utilized to cover broader problem classes (cf. \cite{hahn2019interval, haesaert2021formal}). Further, joint chance constrained policies can yield a hedonistic behaviour, incentively violating safety specifications to exploit low cost trajectories \cite{schmid2023Jccoc}. Different definitions of safety might help overcoming these issues and are subject of future research. 
}



\bibliographystyle{ieeetr}        
\bibliography{root}

\begin{thebibliography}{10}

\bibitem{pitchford2007uncertainty}
J.~W. Pitchford, E.~A. Codling, and D.~Psarra, ``Uncertainty and sustainability
  in fisheries and the benefit of marine protected areas,'' {\em Ecol. Model.},
  vol.~207, no.~2-4, 2007.

\bibitem{fochesato2022data}
M.~Fochesato and J.~Lygeros, ``Data-driven distributionally robust bounds for
  stochastic model predictive control,'' in {\em CDC}, IEEE, 2022.

\bibitem{schmid2022real}
N.~Schmid, J.~Gruner, H.~S. Abbas, and P.~Rostalski, ``{A real-time GP based
  MPC for quadcopters with unknown disturbances},'' in {\em ACC}, IEEE, 2022.

\bibitem{blackmore2011chance}
L.~Blackmore, M.~Ono, and B.~C. Williams, ``Chance-constrained optimal path
  planning with obstacles,'' {\em IEEE T-RO}, vol.~27, no.~6, 2011.

\bibitem{Ono_2}
Y.~K. M.~Ono, M.~Pavone and J.~Balaram, ``Chance-constrained dynamic
  programming with application to risk-aware robotic space exploration,'' {\em
  Auton. Robots}, 2015.

\bibitem{schmid2023Jccoc}
N.~Schmid, M.~Fochesato, S.~H.~Q. Li, T.~Sutter, and J.~Lygeros, ``Computing
  optimal joint chance constrained control policies,'' {\em arXiv preprint
  arXiv:2312.10495}, 2023.

\bibitem{Wang}
K.~Wang and S.~Gros, ``Solving mission-wide chance-constrained optimal control
  using dynamic programming,'' in {\em CDC}, IEEE, 2022.

\bibitem{paulson}
J.~A. Paulson, E.~A. Buehler, R.~D. Braatz, and A.~Mesbah, ``Stochastic model
  predictive control with joint chance constraints,'' {\em IJC}, vol.~93,
  no.~1, 2020.

\bibitem{patil_1}
A.~Patil and T.~Tanaka, ``Upper and lower bounds for end-to-end risks in
  stochastic robot navigation,'' {\em IFAC-PapersOnLine}, vol.~56, no.~2, 2023.

\bibitem{wang2020non}
A.~Wang, A.~Jasour, and B.~C. Williams, ``Non-gaussian chance-constrained
  trajectory planning for autonomous vehicles under agent uncertainty,'' {\em
  IEEE RA-L}, vol.~5, no.~4, 2020.

\bibitem{nemirovski2007convex}
A.~Nemirovski and A.~Shapiro, ``Convex approximations of chance constrained
  programs,'' {\em SIOPT}, vol.~17, no.~4, 2007.

\bibitem{hokayem2013chance}
P.~Hokayem, D.~Chatterjee, and J.~Lygeros, ``Chance-constrained lqg with
  bounded control policies,'' in {\em CDC}, IEEE, 2013.

\bibitem{thorpe_1}
A.~Thorpe, T.~Lew, M.~Oishi, and M.~Pavone, ``Data-driven chance constrained
  control using kernel distribution embeddings,'' in {\em L4DC}, PMLR, 2022.

\bibitem{huang2021risk}
X.~Huang, M.~Feng, A.~Jasour, G.~Rosman, and B.~Williams, ``Risk conditioned
  neural motion planning,'' in {\em IROS}, IEEE, 2021.

\bibitem{haesaert2021formal}
S.~Haesaert, P.~Nilsson, and S.~Soudjani, ``Formal multi-objective synthesis of
  continuous-state {MDPs},'' in {\em ACC}, IEEE, 2021.

\bibitem{hahn2019interval}
E.~M. Hahn, V.~Hashemi, H.~Hermanns, M.~Lahijanian, and A.~Turrini, ``Interval
  {Markov} decision processes with multiple objectives: from robust strategies
  to pareto curves,'' {\em TOMACS}, vol.~29, no.~4, 2019.

\bibitem{de2003linear}
D.~P. De~Farias and B.~Van~Roy, ``The linear programming approach to
  approximate dynamic programming,'' {\em Oper. Res.}, vol.~51, no.~6, 2003.

\bibitem{aumann_1}
R.~J. Aumann, {\em Mixed and Behavior Strategies in Infinite Extensive Games}.
\newblock Princeton University Press, 1964.

\bibitem{abate2008probabilistic}
A.~Abate, M.~Prandini, J.~Lygeros, and S.~Sastry, ``Probabilistic reachability
  and safety for controlled discrete time stochastic hybrid systems,'' {\em
  Automatica}, vol.~44, no.~11, 2008.

\bibitem{bertsekas2012dynamic}
D.~Bertsekas, {\em Dynamic programming and optimal control: Volume I}, vol.~4.
\newblock Athena scientific, 2012.

\bibitem{ioffe2005generic}
A.~Ioffe, R.~Lucchetti, {\em et~al.}, ``Generic well-posedness in minimization
  problems,'' in {\em Abstr. Appl. Anal.}, vol.~2005, Hindawi, 2005.

\bibitem{hernandez2012discrete}
O.~Hern{\'a}ndez-Lerma and J.~B. Lasserre, {\em Discrete-time Markov control
  processes: basic optimality criteria}, vol.~30.
\newblock Springer Science \& Business Media, 2012.

\bibitem{schmid2022probabilistic}
N.~Schmid and J.~Lygeros, ``Probabilistic reachability and invariance
  computation of stochastic systems using linear programming,'' {\em
  IFAC-PapersOnLine}, vol.~56, no.~2, 2023.

\bibitem{kariotoglou2017linear}
N.~Kariotoglou, M.~Kamgarpour, T.~H. Summers, and J.~Lygeros, ``{The linear
  programming approach to reach-avoid problems for Markov decision
  processes},'' {\em JAIR}, vol.~60, 2017.

\end{thebibliography}

\end{document}